\documentclass[12pt,twoside,leqno]{article}
\usepackage{amsmath}
\usepackage{amssymb}
\usepackage{amsxtra}
\usepackage{amscd}
\usepackage{amsthm}
\usepackage[mathscr]{eucal}

\theoremstyle{plain}
\newtheorem{thm}[subsection]{Theorem}

\newtheorem{cor}[subsection]{Corollary}
\newtheorem{lem}[subsection]{Lemma}

\theoremstyle{definition}

\newtheorem{rem}[subsection]{Remark}
\newtheorem{para}[subsection]{}

\newenvironment{pf}{\proof[\proofname]}{\endproof}

\usepackage{fancyhdr}
\usepackage{amssymb}
\usepackage{graphicx}

\begin{document}
\title{Toroidal compactifications and Borel--Serre compactifications}
\author
{Kazuya Kato, Chikara Nakayama, Sampei Usui}

\maketitle

\renewcommand{\mathbb}{\bold}

\newcommand\Cal{\mathcal}
\newcommand\define{\newcommand}

\define\gp{\mathrm{gp}}%
\define\fs{\mathrm{fs}}%
\define\an{\mathrm{an}}%
\define\mult{\mathrm{mult}}%
\define\Ker{\mathrm{Ker}\,}%
\define\Coker{\mathrm{Coker}\,}%
\define\Hom{\mathrm{Hom}\,}%
\define\Ext{\mathrm{Ext}\,}%
\define\rank{\mathrm{rank}\,}%
\define\gr{\mathrm{gr}}%
\define\cHom{\Cal{Hom}}
\define\cExt{\Cal Ext\,}%

\define\cB{\Cal B}
\define\cC{\Cal C}
\define\cD{\Cal D}
\define\cO{\Cal O}
\define\cS{\Cal S}
\define\cM{\Cal M}
\define\cG{\Cal G}
\define\cH{\Cal H}
\define\cE{\Cal E}
\define\cF{\Cal F}
\define\cN{\Cal N}
\define\fF{\frak F}
\define\Dc{\check{D}}
\define\Ec{\check{E}}

\newcommand{\N}{{\mathbb{N}}}
\newcommand{\Q}{{\mathbb{Q}}}
\newcommand{\Z}{{\mathbb{Z}}}
\newcommand{\R}{{\mathbb{R}}}
\newcommand{\C}{{\mathbb{C}}}
\newcommand{\bN}{{\mathbb{N}}}
\newcommand{\bQ}{{\mathbb{Q}}}
\newcommand{\bF}{{\mathbb{F}}}
\newcommand{\bZ}{{\mathbb{Z}}}
\newcommand{\bP}{{\mathbb{P}}}
\newcommand{\bR}{{\mathbb{R}}}
\newcommand{\bC}{{\mathbb{C}}}
\newcommand{\bbQ}{{\bar \mathbb{Q}}}
\newcommand{\ol}[1]{\overline{#1}}
\newcommand{\too}{\longrightarrow}
\newcommand{\respect}{\rightsquigarrow}
\newcommand{\compatible}{\leftrightsquigarrow}
\newcommand{\upc}[1]{\overset {\lower 0.3ex \hbox{${\;}_{\circ}$}}{#1}}
\newcommand{\Gmlog}{\bG_{m, \log}}
\newcommand{\Gm}{\bG_m}
\newcommand{\ep}{\varepsilon}
\newcommand{\Spec}{\operatorname{Spec}}
\newcommand{\val}{{\mathrm{val}}} 
\newcommand{\n}{\operatorname{naive}}
\newcommand{\bs}{\operatorname{\backslash}}
\newcommand{\Gal}{\operatorname{{Gal}}}
\newcommand{\gal}{{\rm {Gal}}({\bar \Q}/{\Q})}
\newcommand{\galp}{{\rm {Gal}}({\bar \Q}_p/{\Q}_p)}
\newcommand{\gall}{{\rm{Gal}}({\bar \Q}_\ell/\Q_\ell)}
\newcommand{\wep}{W({\bar \Q}_p/\Q_p)}
\newcommand{\wel}{W({\bar \Q}_\ell/\Q_\ell)}
\newcommand{\Ad}{{\rm{Ad}}}
\newcommand{\BS}{{\rm {BS}}}
\newcommand{\even}{\operatorname{even}}
\newcommand{\End}{{\rm {End}}}
\newcommand{\odd}{\operatorname{odd}}
\newcommand{\GL}{\operatorname{GL}}
\newcommand{\np}{\text{non-$p$}}
\newcommand{\g}{{\gamma}}
\newcommand{\G}{{\Gamma}}
\newcommand{\Lam}{{\Lambda}}
\newcommand{\La}{{\Lambda}}
\newcommand{\lam}{{\lambda}}
\newcommand{\la}{{\lambda}}
\newcommand{\uL}{{{\hat {L}}^{\rm {ur}}}}
\newcommand{\uQp}{{{\hat \Q}_p}^{\text{ur}}}
\newcommand{\sel}{\operatorname{Sel}}
\newcommand{\dt}{{\rm{Det}}}
\newcommand{\Sig}{\Sigma}
\newcommand{\fil}{{\rm{fil}}}
\newcommand{\SL}{{\rm{SL}}}
\newcommand{\spl}{{\rm{spl}}}
\newcommand{\st}{{\rm{st}}}
\newcommand{\Isom}{{\rm {Isom}}}
\newcommand{\Mor}{{\rm {Mor}}}
\newcommand{\bg}{\bar{g}}
\newcommand{\id}{{\rm {id}}}
\newcommand{\cone}{{\rm {cone}}}
\newcommand{\al}{a}
\newcommand{\ChL}{{\cal{C}}(\La)}
\newcommand{\Image}{{\rm {Image}}}
\newcommand{\toric}{{\operatorname{toric}}}
\newcommand{\torus}{{\operatorname{torus}}}
\newcommand{\Aut}{{\rm {Aut}}}
\newcommand{\Qp}{{\mathbb{Q}}_p}
\newcommand{\barQp}{{\mathbb{Q}}_p}
\newcommand{\Qpur}{{\mathbb{Q}}_p^{\rm {ur}}}
\newcommand{\Zp}{{\mathbb{Z}}_p}
\newcommand{\Zl}{{\mathbb{Z}}_l}
\newcommand{\Ql}{{\mathbb{Q}}_l}
\newcommand{\Qlur}{{\mathbb{Q}}_l^{\rm {ur}}}
\newcommand{\F}{{\mathbb{F}}}
\newcommand{\eps}{{\epsilon}}
\newcommand{\epsLa}{{\epsilon}_{\La}}
\newcommand{\epsLaVxi}{{\epsilon}_{\La}(V, \xi)}
\newcommand{\epsOLaVxi}{{\epsilon}_{0,\La}(V, \xi)}
\newcommand{\Qplin}{{\mathbb{Q}}_p(\mu_{l^{\infty}})}
\newcommand{\otimesQplin}{\otimes_{\Qp}{\mathbb{Q}}_p(\mu_{l^{\infty}})}
\newcommand{\galFl}{{\rm{Gal}}({\bar {\Bbb F}}_\ell/{\Bbb F}_\ell)}
\newcommand{\gallur}{{\rm{Gal}}({\bar \Q}_\ell/\Q_\ell^{\rm {ur}})}
\newcommand{\galFF}{{\rm {Gal}}(F_{\infty}/F)}
\newcommand{\galFv}{{\rm {Gal}}(\bar{F}_v/F_v)}
\newcommand{\galF}{{\rm {Gal}}(\bar{F}/F)}
\newcommand{\epsVxi}{{\epsilon}(V, \xi)}
\newcommand{\epsOVxi}{{\epsilon}_0(V, \xi)}
\newcommand{\plim}{\lim_
{\scriptstyle 
\longleftarrow \atop \scriptstyle n}}
\newcommand{\sig}{{\sigma}}
\newcommand{\ga}{{\gamma}}
\newcommand{\del}{{\delta}}
\newcommand{\Vss}{V^{\rm {ss}}}
\newcommand{\Bst}{B_{\rm {st}}}
\newcommand{\Dpst}{D_{\rm {pst}}}
\newcommand{\Dcrys}{D_{\rm {crys}}}
\newcommand{\DdR}{D_{\rm {dR}}}
\newcommand{\Fin}{F_{\infty}}
\newcommand{\Kla}{K_{\lambda}}
\newcommand{\Ola}{O_{\lambda}}
\newcommand{\Mla}{M_{\lambda}}
\newcommand{\Det}{{\rm{Det}}}
\newcommand{\Sym}{{\rm{Sym}}}
\newcommand{\LaSa}{{\La_{S^*}}}
\newcommand{\cX}{{\cal {X}}}
\newcommand{\MHG}{{\frak {M}}_H(G)}
\newcommand{\tauMla}{\tau(M_{\lambda})}
\newcommand{\Fvur}{{F_v^{\rm {ur}}}}
\newcommand{\Lie}{{\rm {Lie}}}
\newcommand{\cL}{{\cal {L}}}
\newcommand{\cW}{{\cal {W}}}
\newcommand{\fq}{{\frak {q}}}
\newcommand{\cont}{{\rm {cont}}}
\newcommand{\SC}{{SC}}
\newcommand{\Om}{{\Omega}}
\newcommand{\dR}{{\rm {dR}}}
\newcommand{\crys}{{\rm {crys}}}
\newcommand{\hatSig}{{\hat{\Sigma}}}
\newcommand{\rdet}{{{\rm {det}}}}
\newcommand{\ord}{{{\rm {ord}}}}
\newcommand{\BdR}{{B_{\rm {dR}}}}
\newcommand{\BdRO}{{B^0_{\rm {dR}}}}
\newcommand{\Bcrys}{{B_{\rm {crys}}}}
\newcommand{\Qw}{{\mathbb{Q}}_w}
\newcommand{\barkappa}{{\bar{\kappa}}}
\newcommand{\cP}{{\Cal {P}}}
\newcommand{\cZ}{{\Cal {Z}}}
\newcommand{\oppLa}{{\Lambda^{\circ}}}
\newcommand{\add}{{{\rm {add}}}}
\newcommand{\red}{{{\rm {red}}}}
\newcommand\twoheaddownarrow{\mathrel{\rotatebox[origin=c]{90}{$\twoheadleftarrow$}}}

\begin{abstract} 
  We discuss connections of toroidal compactifications and Borel--Serre compactifications in view of the fundamental diagram of extended period domains.  We give a complement to a work of Goresky--Tai. 
\end{abstract}

\renewcommand{\thefootnote}{\fnsymbol{footnote}}
\footnote[0]{MSC2020: Primary 14A21; Secondary 14D07, 32G20} 


\footnote[0]{Keywords: Hodge theory, moduli, %
toroidal compactification, Borel--Serre compactification}

\setcounter{section}{-1}
\section{Introduction}\label{s:intro}
  In this paper,
as an application of the fundamental diagram of extended period domains introduced in our series \cite{KNU1}--\cite{KNU5} of papers, we give a complement to the work of Goresky and Tai \cite{GT} on the relation of the toroidal compactification and the reductive Borel--Serre compactification.

  Let $G$ be a reductive algebraic group over $\bQ$ and let $h_0: S_{\bC/\bR} \to G_{\bR}$ be a homomorphism satisfying certain conditions which are necessary to consider the period domain $D$ associated to $(G, h_0)$.  
  Then the main part of the associated  fundamental diagram in \cite{KNU5} is the following. 
$$\begin{matrix}
\Gamma \bs D_{\Sig, [\val]} &\twoheadleftarrow & D^{\sharp}_{\Sig, [\val]} & \overset{\psi}\to & D_{\SL(2),\val} & \overset{\eta}\to & D_{\BS,\val}\\
\twoheaddownarrow && \twoheaddownarrow && \twoheaddownarrow && \twoheaddownarrow\\
\Gamma \bs D_{\Sig, [:]}&\twoheadleftarrow& D^{\sharp}_{\Sig,[:]} & \overset{\psi}\to & D_{\SL(2)} && D_{\BS}\\
\twoheaddownarrow && 
\twoheaddownarrow  
&&&&\twoheaddownarrow \\
\Gamma \bs D_{\Sig} &\twoheadleftarrow & D^{\sharp}_{\Sig} &&&& D_{\BS}^{\mathrm{red}}
\end{matrix}$$
  Here the notation is as in {\rm \cite{KNU5}} Section 4.
  In particular, $\Sig$ is a weak fan and $\Gamma$ is a neat semi-arithmetic subgroup (\ref{Gamma}) of $G'(\bQ)$ which is strongly compatible with $\Sig$. 
  ($G'$ denotes the commutator group of $G$.) 
  
  In this diagram, all maps are continuous. All vertical maps and the maps to the left direction are surjective and the targets of these maps have the quotient topologies. The vertical arrows  except the map $D_{\BS}\to D^{\red}_{\BS}$ are proper. 

  Assume further that $(G,h_0)$ comes from some Shimura data and that $\Gamma$ is an arithmetic subgroup (\ref{Gamma}). 
  Then we have the following two kinds of compactifications of the Shimura variety $\Gamma \bs D$. 
   On one hand, there is a $\Sig$ such that $\Gamma\bs D_{\Sig}$ is compact and this is a toroidal compactification of $\Gamma \bs D$. 
   On the other hand, we have the  reductive Borel--Serre compactification  $\Gamma \bs D_{\BS}^{\mathrm{red}}$ of $\Gamma \bs D$, which is a quotient of the  
   Borel--Serre compactification $\Gamma \bs D_{\BS}$ of $\Gamma \bs D$.
   
   Now, Goresky and Tai constructed a space which connects, after replacing $\Sig$ by a finer subdivision, these two compactifications 
of $\Gamma \bs D$ in the following fashion.   
  That is, they show that, after replacing $\Sig$ by a sufficiently finer subdivision $\Sig'$, there are a compact topological space $T$ which contains $\Gamma \bs D$ as a dense open subspace, and continuous surjective maps $f: T\to \Gamma\bs D_{\Sig'}$ and  $g: T\to \Gamma \bs D^{\mathrm{red}}_{\BS}$ such that $f$ is a homotopy equivalence and such that $f$ and $g$ induce the identity map of $\Gamma \bs D$: 
$$
\begin{CD}
\Gamma\bs D_{\Sig'} @<\mathrm{hom.\ eq.}<< T @>>> \Gamma \bs D^{\mathrm{red}}_{\BS}.
\end{CD}
$$
  From this, the existence of the canonical maps 
$$H^m(\Gamma \bs D^{\mathrm{red}}_{\BS},A) \to H^m(\Gamma \bs D_{\Sig}, A)$$ 
with $A$ being an abelian group for any $m$, $A$ being a group for $m=0,1$, and $A$ being a set for $m=0$ follows. 

  In this paper, we have the following complement to the work of Goresky and Tai (cf.\ \ref{app1}). 
  As we have shown in \cite{KNU5} 3.8, under the above assumption, there is a unique continuous map $D_{\SL(2)}\to D_{\BS}$ which makes the above fundamental diagram still commutative, and the composition $D^{\sharp}_{\Sig, [:]}\overset {\psi} \to D_{\SL(2)}\to D_{\BS}$ induces the surjection of the quotient spaces $\Gamma \bs D_{\Sig, [:]}\to \Gamma \bs D^{\mathrm{red}}_{\BS}$.
  On the other hand, as we will prove in this paper, 
the map $\Gamma \bs D_{\Sig, [:]} \to \Gamma \bs D_{\Sig}$ in the fundamental diagram is proper, surjective, and a weak homotopy equivalence (Theorem \ref{homoto2} (3)): 
$$
\begin{CD}
\Gamma\bs D_{\Sig} @<\mathrm{weak\ hom.\  eq.}<< \Gamma \bs D_{\Sig, [:]}@>>>\Gamma \bs D^{\mathrm{red}}_{\BS}.
\end{CD}
$$
  That is, compared to \cite{GT}, 
we do not need a subdivision of $\Sig$ and our intermediate space is a concrete one, that is, $\Gamma \bs D_{\Sig, [:]}$. 
  Thus we can give an alternative construction of the above canonical maps between cohomologies (Corollary \ref{c:coh}). 
  (Probably 
$\Gamma \bs D_{\Sig, [:]} \to \Gamma \bs D_{\Sig}$ is in fact a homotopy equivalence.  
For this, it is enough to see that $\Gamma \bs D_{\Sig, [:]}$ is triangulizable, which is plausible, though we don't try to show it.)

  Further, even not necessarily in the situation of Shimura data, we have 
$$
\begin{CD}
\Gamma\bs D_{\Sig} @<\mathrm{same\ cohomologies}<< \Gamma \bs D_{\Sig, [\val]}@>>>\Gamma \bs D^{\mathrm{red}}_{\BS},
\end{CD}
$$
where the first map induces the isomorphisms of cohomologies including $H^1$ with nonabelian coefficients, and we still have the above canonical maps between cohomologies. 

  The organization of this paper is as follows. 
  Section \ref{s:rev} is a preparation.
  In Section \ref{s:con}, we prove the basic result that the fibers of the natural projection from the ratio space are weakly contractible.
  Based on this result, the complement for Goresky--Tai alluded in the above is an easy consequence of the existence of our fundamental diagram, which is explained in Section \ref{s:GT}. 
  Section \ref{s:val} is a valuative version of the results in Section \ref{s:con}.
  Section \ref{s:general} is a generalization of the results in Section \ref{s:GT} to the case where $(G,h_0)$ does not necessarily come from Shimura data.

\medskip

  K.\ Kato was 
partially supported by NFS grants DMS 1303421, DMS 1601861, and DMS 2001182.
C.\ Nakayama was 
partially supported by JSPS Grants-in-Aid for Scientific Research (B) 23340008, (C) 16K05093, and (C) 21K03199.
S.\ Usui was 
partially supported by JSPS Grants-in-Aid for Scientific Research (B) 23340008 and (C) 17K05200.

\bigskip

\section{Preliminaries} 
\label{s:rev}
  In this section, we gather terminology and some well-known facts used in the following sections. 

\begin{para}
\label{coh}
  In this paper, we say that a continuous map $f: X \to Y$ of topological spaces {\it induces the isomorphisms of cohomologies} 
if for any locally constant sheaf $F$ of sets (resp.\ groups, resp.\ abelian groups) on $Y$, the induced 
map $H^m(Y, F) \to H^m(X, f^{-1}F)$ is an isomorphism for $m=0$ (resp.\ $m=0,1$, resp.\ $m\in \bZ$).
  Note that we include the nonabelian coefficients here, which is essential in the following. 
\end{para}

\begin{para}
  A topological space is said to be {\it cohomologically contractible} 
if the map to a point induces the isomorphisms of cohomologies in the above sense, i.e., its cohomology groups with abelian and nonabelian coefficients are all trivial. 

  A topological space is said to be {\it weakly contractible} 
if the map to a point is a weak homotopy equivalence, i.e., its $\pi_n$ ($n\geq 0$)  are all trivial. 

  A topological space is said to be {\it locally contractible} 
if each point has a fundamental system of contractible open 
neighborhoods. 

  A continuous map $f: X \to Y$ of topological spaces is said to be {\it proper} if it is universally closed and separated. 

  A topological space is said to be {\it compact} if the map to a point is proper. 
  Hence a compact space is Hausdorff. 
\end{para}

\begin{para}
\label{Hure}
  By the Hurewicz theorem, 
a continuous map $f$ of locally contractible spaces 
induces the isomorphisms of cohomologies if and only if it 
is a weak homotopy equivalence (cf.\ \cite{Q} Chapter II, Section 3, Proposition 4).

  In particular, a locally contractible space is cohomologically contractible if and only if it is weakly contractible. 
\end{para}

\begin{para}
\label{pbc}
  A proper map $f: X \to Y$ induces the isomorphisms of cohomologies if each fiber is cohomologically contractible. 
  In fact, let $F$ be as in \ref{coh}. 
  By the proper base change theorem (\cite{KU} Theorem A2.1), 
$R^qf_*G$, where $G=f^{-1}F$ coincides with $F$ if $q=0$ and vanishes if otherwise. 
  Hence $f$ induces the isomorphisms of cohomologies by Leray spectral sequences and the exact sequence of pointed sets 
$H^1(Y, f_*G) \to H^1(X, G) \to H^0(Y, R^1f_*G)$, where the first map is injective. 
\end{para}

\begin{para}
\label{Gamma} 
Let $G$ be a linear algebraic group over $\Q$ and let $\Gamma$ be a subgroup of $G(\Q)$. Then $\Gamma$ is called an {\it arithmetic subgroup} of $G(\Q)$ if there are $n\geq 1$ and  an injective homomorphism $\rho: G\to\GL(n)$ such that $\Gamma$ is a subgroup of $\{g \in G(\Q)\;|\; \rho(g)\in \GL(n,\Z)\}$ of finite index. As in \cite{KNU5}, 
  we call $\Gamma$ a {\it semi-arithmetic subgroup} of $G(\Q)$ if there are $n\geq 1$ and  an injective homomorphism $\rho: G\to\GL(n)$ such that $\rho(\Gamma)\subset \GL(n,\Z)$.
\end{para}

\section{Weak contractibility}
\label{s:con}
  In this section, we prove that the projection from the space of ratios is a weak homotopy equivalence.

\begin{para}
  Let
$S$ be a locally ringed space over $E=\R$ or $\C$ with an fs log structure $M_S$ satisfying the following conditions (i) and (ii).

(i) For every $s\in S$, the natural homomorphism 
$E\to \cO_{S,s}/m_s$ is an isomorphism, where $m_s$ is the maximal ideal of $\cO_{S,s}$.

(ii) For every open set $U$ of $S$ and for every $f\in \cO_S(U)$, the map $U\to E\;;\; s \mapsto f(s)$ is continuous. Here $f(s)$ is the image of $f$ in $\cO_{S,s}/m_s=E$. 

Then we have a topological space $S_{[:]}$ over $S$, called  the {\it space of ratios}.  
  See \cite{KNU4} 4.2 and \cite{KNU5} 4.3. 
\end{para}

In this section we will prove the following theorem, which is the basic result in this paper.

\begin{thm}\label{homoto1}
  Let $S$ be as above. 
  Then the following hold{\rm :}

$(1)$ All fibers of $S_{[:]}\to S$ are weakly contractible and locally contractible. 

$(2)$ The map $S_{[:]}\to S$ induces the isomorphisms of cohomologies ($\ref{coh}$). 
\end{thm}

  Since $S_{[:]} \to S$ is proper, by \ref{Hure} and \ref{pbc}, (2) is reduced to (1). 

(1) follows from the following theorem whose proof occupies the rest of this section. 
Note that the fiber of $S_{[:]}\to S$ on $s\in S$  is identified with the space of ratios $R(\cS)$ of the sharp fs monoid $\cS=(M_S/\cO_S^\times)_s$, which is defined in \cite{KNU4} 4.1.

\begin{thm}\label{contra1}  For a sharp fs monoid $\cS$, the space of ratios $R(\cS)$ is weakly contractible and locally contractible.
\end{thm}

\begin{para} Let $\Phi=\{\cS^{(i)}\; |\; 0\leq i\leq n\}$, where $\cS^{(i)}$ are faces of $\cS$ such that  $\cS=\cS^{(0)}\supsetneq \cS^{(1)}\supsetneq \dots\supsetneq \cS^{(n)}=\{1\}$. Then by  \cite{KNU4} Corollary 4.2.17,
the $\Phi$-part   $R(\cS)(\Phi)$  of $R(\cS)$, which is an open set of $R(\cS)$, has the following description as  a topological space.   For $0\leq i\leq n-1$, take $a_i\in \cS^{(i)}$ such that $a_i\notin \cS^{(i+1)}$. Then  $R(\cS)(\Phi)$ is homeomorphic to the subspace of $\prod_{i=0}^{n-1}\Hom(\cS^{(i)}, \R_{\geq 0}^{\add})$ consisting of all elements $N=(N_i)_{0\leq i \leq n-1}$ satisfying  the following conditions (i), (ii), and (iii).

(i) $N_i(a_i)=1$ for $0\leq i\leq n-1$.

(ii) For $0\leq i\leq n-1$, the kernel of $N_i$ coincides with $\cS^{(j)}$ for some $j$ such that $i< j\leq n$. 

(iii) If $0\leq i < j\leq n-1$, the restriction of $N_i$ to $\cS^{(j)}$ coincides with $N_i(a_j)N_j$. 

We will identify $R(\cS)(\Phi)$ with the set of these $N$. 

\end{para}

\begin{para}

Let $R(\cS)(\Phi)_{>0}$ be the subset of $R(\cS)(\Phi)$ consisting of all $(N_i)_{0\leq i\leq n-1}$ such that $\text{Ker}(N_i)=0$ for all $i$ (this is equivalent to $\text{Ker}(N_0)=0$).

Let $H$ be the set of all homomorphisms 
 $h: \cS\to \R_{\geq 0}^{\add}$ such that $\text{Ker}(h)=0$. We have a surjection
 $$\pi: H\to R(\cS)(\Phi)_{>0}$$ given by $h\mapsto N$, where $N_i$ is the restriction of $h(a_i)^{-1}h$ to $\cS^{(i)}$. 
For $h\in H$, we have $h= h(a_0)\pi(h)_0$. Hence  
 $\pi$ induces a homeomorphism $\{h\in H\;|\; h(a_0)=1\}\to R(\cS)(\Phi)_{>0}$, whose inverse map is given by $N\mapsto N_0$.

\end{para}

\begin{para} For $0\leq i\leq n-1$, fix a homomorphism $p_i: \cS\to \cS^{(i)}_{\Q_{\ge_0}}$ whose restriction to $\cS^{(i)}$ is the inclusion map $\cS^{(i)} \to \cS^{(i)}_{\bQ_{\ge_0}}$. 
  The existence of such $p_i$ is seen by the next lemma, which is probably well-known.
             
  For an $N\in R(\cS)(\Phi)$ and a real number $t>0$, define $$\theta_t(N):=\sum_{i=0}^{n-1}   t^i N_i \circ p_i\in H.$$ 
\end{para}

\begin{lem}
\label{l:section}
Let $\cS$ be a sharp fs monoid and let $\cS'$ be a face of $\cS$. 
  Then there is a homomorphism $\cS\to \cS'_{\Q_{\ge0}}$ whose restriction to $\cS'$ is the inclusion map $\cS'\to \cS'_{\Q_{\ge0}}$. 
\end{lem}

\begin{pf} 
  By the induction on $\dim \cS'$, we may assume that $\dim \cS'=\dim \cS-1$. 
  Take a direct decomposition $\cS_{\bQ} = \cS'_{\bQ} \oplus \bQ$ with the first and the second projection being $p$ and $l$ respectively such that $l(\cS) \subset \bQ_{\ge 0}$.
  Let $a \in \cS'$ and consider the homomorphism
$p_a: \cS \to \cS'_\Q \;;\; x\mapsto  p(x) + l(x)a$ (the semi-group law is written additively here).
Then $p_a(x)=x$ for $x\in \cS'$. Hence it is sufficient to prove that there is an $a$ such that $p_a(\cS)\subset \cS'_{\Q_{\ge0}}$.
  
  Let $x_1, \dots, x_n$ be a set of generators of $\cS$.
  If $x_j\in \cS'$, $p_a(x_j)=x_j \in \cS'$.
  If $x_j \notin \cS'$,  $l(x_j) >0$ so that if $a$ is an inner point of $\cS'$, $l(x_j)a$ is an inner point of $\cS'_{\Q_{\ge0}}$. 
  Replacing $a$ by $ca$ for $c\gg 0$, we have $p_a(x_j)=p(x_j)+l(x_j)a\in \cS'_{\Q_{\geq 0}}$ for all $j$.
\end{pf}

\begin{lem}
\label{pfcontra1} 
$R(\cS)(\Phi)$ is contractible. 
\end{lem}

\begin{pf} Fix an $L\in H$. Consider the continuous map   $$f: R(\cS)(\Phi)\times (0, 1] \to R(\cS)(\Phi)\;\;;\; \;
f(N, t)= \pi(t^nL+ (1-t)\theta_t(N)).$$ Then $f(N, 1)= \pi(L)$. It is sufficient to prove that $f$ extends to a continuous map $f: R(\cS)(\Phi)\times [0, 1] \to R(\cS)(\Phi)$ such that $f(N, 0)=N$. Assume $M\in R(\cS)(\Phi)$ converges to $N\in R(\cS)(\Phi)$ and  $t\in (0,1)$ converges to $0$. We prove that $f(M, t)$ converges to $N$.
Let $0\leq i\leq n-1$. Then 
the restriction of $t^nL+(1-t)\theta_t(M)$ to $\cS^{(i)}$ coincides with $t^i(b(t)M_i + E(t))$, where 
$b(t)= (1-t)\sum_{k=0}^i t^{k-i}M_k(a_i)\geq (1-t)M_i(a_i)=1-t$ and $E(t)$ is the restriction of $t^{n-i}L+ (1-t)\sum_{k=i+1}^{n-1} t^{k-i}M_k \circ p_k$ to $\cS^{(i)}$. Hence 
$\pi(t^nL+(1-t)\theta_t(M))_i=  (M_i+b(t)^{-1}E(t))(1+ b(t)^{-1}E(t)(a_i))^{-1}$. This converges to $N_i$ since $b(t) \geq 1-t$  and since $M_i$ and $E(t)$ converge to $N_i$ and $0$, respectively. 

\end{pf}

\begin{para}\label{pfcontra2} 
  We prove that $R(\cS)$ is weakly and locally contractible. 
  First, by \cite{KNU5} Proposition 4.2.19, $R(\cS)$ is a topological manifold with boundaries. 
  In particular, it is locally contractible. 
  Next, $R(\cS)$ has an open covering by $\{R(\cS)(\Phi)\}_{\Phi}$, and for any $\Phi_1$ and $\Phi_2$, we have $R(\cS)(\Phi_1)\cap R(\cS)(\Phi_2)= R(\cS)(\Phi_1 \cap \Phi_2)$.
  Together with this, Lemma \ref{pfcontra1} implies that $R(\cS)$ is cohomologically contractible.
  Hence, by \ref{Hure}, $R(\cS)$ is weakly contractible. 
\end{para}

\begin{rem}
  It is plausible that $R(\cS)$ is in fact contractible and even is a homeomorphic to a closed ball.
  At least, $R(\cS)$ is a compact Hausdorff topological manifold with boundaries whose interior is homeomorphic to an open ball.  
  But the authors do not know if this implies that $R(\cS)$ is homeomorphic to a closed ball.
\end{rem}

\section{On the work of Goresky and Tai}
\label{s:GT}

\begin{para}
\label{setting}
  Let $G$ be a reductive algebraic group over $\bQ$. 
  Let $h_0: S_{\bC/\bR} \to G_{\bR}$ be a homomorphism as in \cite{KNU5} 1.2.13. 
  Assume that $h_0$ is $\bR$-polarizable (\cite{KNU5} 1.5.2). 
  Let $D=D(G,h_0)$ be the period domain associated to $(G,h_0)$. 
  Let $\Sig$ be a weak fan in $\Lie(G')$ (\cite{KNU5} 4.1.7), where $G'$ denotes the commutator group of $G$.  
  Let $\Gamma$ be a neat semi-arithmetic subgroup (\ref{Gamma}) of $G'(\Q)$.
  Assume that $\Sig$ and $\Gamma$ are strongly compatible. 

Below, the {\it Shimura data case} 
(\cite{D} 1.5)
means the case where the Hodge type of $\Lie(G_{\bR})$ via $h_0$ is in $\{(1,-1), (0, 0), (-1,1)\}$ (as in Shimura data).  
\end{para}

\begin{thm}\label{homoto2}
  Let the notation and the assumptions be as above.
  Then the following hold{\rm :}

$(1)$ All fibers of $\Gamma \bs D_{\Sig,[:]}\to \Gamma\bs D_{\Sig}$ are weakly contractible and locally contractible. 

$(2)$ The map $\Gamma \bs D_{\Sig,[:]}\to \Gamma\bs D_{\Sig}$ induces the isomorphisms of cohomologies ($\ref{coh}$). 

$(3)$ In the Shimura data case ($\ref{setting}$), 
$\Gamma \bs D_{\Sig,[:]}\to \Gamma\bs D_{\Sig}$ is also a weak homotopy equivalence.
\end{thm}

\begin{pf}
  (1) is by Theorem \ref{homoto1} (1). 
  (2) follows from (1) by \ref{Hure} and \ref{pbc}. 
  Then, by \ref{Hure}, to see (3), it is enough to show that, in the Shimura data case, $\Gamma \bs D_{\Sig,[:]}$ is locally triangulizable so that locally contractible. 
  Take a local chart of the fs log complex analytic space $\Gamma \bs D_{\Sig}$. 
  Then, by \cite{KNU4} Proposition 4.2.14, $\Gamma \bs D_{\Sig,[:]}$ is locally a fiber product of two analytic maps over a real manifold with singular corners.   Hence it is locally a semianalytic subset of an Euclidean space so that locally triangulizable by a theorem of Hironaka.
\end{pf}

\begin{para}\label{app1}
  As is explained in Introduction, 
this Theorem \ref{homoto2} gives a new interpretation to the main result 
of the work of Goresky and Tai on the relation between the toroidal compactification and the reductive Borel--Serre compactification, which we will explain.

  Let the reductive Borel--Serre space $D_{\BS}^{\mathrm{red}}$ be the quotient of $D_{\BS}$ by the following equivalence relation. 
For $p_1=(P_1, Z_1), p_2=(P_2, Z_2)\in D_{\BS}$, $p_1\sim p_2$ if and only if $P_1=P_2$ and $P_{1,u}Z_1=P_{2,u}Z_2$, where $P_{i,u}$ denotes the unipotent radical of the parabolic subgroup $P_i$ $(i=1,2)$. 

  In the case where $D$ is a Griffiths domain \cite{G}, this $D^{\red}_{\BS}$ is denoted by $D^{\flat}_{\BS}$ in \cite{KU} Section 9.
  If $\Gamma$ is semi-arithmetic (resp.\ arithmetic), then $\Gamma \bs D^{\red}_{\BS}$ is Hausdorff (resp.\ compact). 
  The proof is similar to that in \cite{KU} Section 9.2. 

  In the Shimura data case,  for an arithmetic $\Gamma$, 
$\Gamma\bs D^{\red}_{\BS}$ is the well-known reductive Borel--Serre compactification of the symmetric domain $\Gamma \bs D$ studied in \cite{Z}.

  Now assume that we are in the Shimura data case. 
  Then we have the natural morphism $D_{\SL(2)}\to D_{\BS}$ in the fundamental diagram (\cite{KNU5} Theorem 3.8.2) and hence the CKS map $D_{\Sig,[:]}^{\sharp}\to D_{\SL(2)}$ induces a continuous map $D_{\Sig, [:]}^{\sharp}\to D_{\BS}$, 
which 
induces a continuous map 
$\Gamma \bs D_{\Sig,[:]}\to \Gamma\bs D^{\red}_{\BS}$.

  On the other hand, the map $\Gamma \bs D_{\Sig,[:]}\to \Gamma \bs D_{\Sig}$ is 
a weak homotopy equivalence (Theorem \ref{homoto2} (3)).
  Thus the formation of these maps are regarded as a variant of the construction by Goresky and Tai. 
\end{para}

\begin{cor}
\label{c:coh}
  Let $A$ be a set (resp.\ a group, resp.\ an abelian group).
  Then, in the Shimura data case ($\ref{setting}$), there is a canonical map 
$$H^m(\Gamma \bs D^{\mathrm{red}}_{\BS},A) \to H^m(\Gamma \bs D_{\Sig}, A)$$ 
for $m=0$ (resp.\ $m=0,1$, resp.\ $m\in \bZ$).
\end{cor}

\section{Valuative version}
\label{s:val}

  In this section, we consider the valuative version of the results in Section \ref{s:con}. 

\begin{thm}\label{homoto3} 
  Let $S$ be either an object of $\cB'_\R(\log)$  ({\rm \cite{KNU2} 3.1, \cite{KNU4} 1.3, \cite{KNU5} 3.4.1}) or an fs log topological space whose structural sheaf of rings is the sheaf of all continuous $\R$-valued functions. Let $S_{\val}$ (resp. $S_{[\val]}$) be the topological space defined as in {\rm \cite{KNU4} 3.1.2} and {\rm \cite{KNU5} 4.4} 
(resp. {\rm \cite{KNU4} 4.3} and {\rm \cite{KNU5} 4.7}). 
Then the following hold{\rm:}

$(1)$  All fibers of $S_{\val}\to S$ and all fibers of $S_{[\val]}\to S$ are cohomologically contractible. 

$(2)$ 
  The map $S_{\val}\to S$ and the map $S_{[\val]}\to S$ induce the isomorphisms of cohomologies ($\ref{coh}$). 
\end{thm}

\begin{pf} 
  Since $S_{\val}\to S$ and $S_{[\val]} \to S$ are proper, (2) is reduced to (1) by \ref{pbc}. 
  Because $S_{[\val]}\to S$ is the composition $S_{[\val]}\to S_{[:]}\to S$ and $S_{[\val]}=(S_{[:]})_{\val}$ for the new log structure of $S_{[:]}$, the $S_{[\val]}$ case is reduced to the $S_{\val}$ case by Theorem \ref{homoto1} and the proper base change theorem.

  Thus it suffices to show  the $S_{\val}$ case of (1). 
  It is sufficient to prove that for $f,g\in \cS$, the log blowing-up $B$ of $S$ by $(f,g)$ is contractible, 
where $\cS$ is a sharp fs monoid and $S$ is a point with the ring $\R$ and the log structure defined by the local homomorphism $\cS\to \R_{\geq 0}$ (cf.\ \cite{KjNk} p.311). 
If either $f/g\in \cS$ or $g/f\in \cS$, $B=S$ and hence $B$ is contractible. Otherwise, $B$ is homeomorphic to the interval $[0, \infty]$ and hence $B$  is contractible. 
\end{pf}

\begin{rem}
  The $S_{\val}$ part of (1) gives an alternative proof of Theorem \ref{homoto1} without constructing homotopies, which we sketch here. 
  It is enough to show that $S_{\val} \to S_{[:]}$ induces the isomorphisms of cohomologies, where $S$ is as in the last paragraph in the proof of the theorem.
  By \ref{pbc}, it is enough to show that any fiber $X$ of the last map is cohomologically contractible. 
  This $X$ is the closure of a stratum $Y$ of $S_{\val}$. 
  For each log modification $S'$ of $S$, the closure of the image of $Y$ in $S'$ is homeomorphic to a closed ball, and $X$ is the inverse limit of these closures.
  Hence $X$ is cohomologically contractible. 
\end{rem}

The above proof of Theorem \ref{homoto3} shows the following. 

\begin{thm}
\label{t:realbu}
For an $S$ as in Theorem $\ref{homoto3}$, for every log modification of $S$, all fibers are cohomologically contractible.
\end{thm}

\begin{rem}
  The fibers in Theorem \ref{t:realbu} are in fact contractible because they are triangulizable by the theory of toric varieties. 
  See \cite{bc}.
  One can ask if it is homeomorphic to a closed ball. 
\end{rem}

\section{Relation to extended period domains}
\label{s:general}

We can generalize \ref{app1} to the non-classical (without Shimura data) situation, where we may not have $D_{\SL(2)} \to D_{\BS}$.
  Here we only assume that $G$ is reductive. 

\begin{para}

  By Theorem \ref{homoto3} (1), all fibers of 
$$\Gamma \bs D_{\Sig, [\val]}\to \Gamma \bs D_{\Sig,[:]}\to \Gamma \bs D_{\Sig}$$ are cohomologically contractible. 
  Hence, by \ref{pbc}, this map induces the isomorphisms of cohomologies ($\ref{coh}$). 
\end{para}

\begin{para}
On the other hand, we have continuous maps 
$$\Gamma \bs D_{\Sig, [\val]}\to 
\Gamma \bs D^{\red}_{\BS}$$
induced by the continuous map 
$$D_{\Sig,[\val]}^{\sharp} \to D_{\SL(2),\val}\to D_{\BS, \val}\to D_{\BS}
\to D^{\red}_{\BS}
$$
in the fundamental diagram.
\end{para}

\begin{cor}
  Let $A$ be a set (resp.\ a group, resp.\ an abelian group).
  Then, there is a canonical map 
$$H^m(\Gamma \bs D^{\mathrm{red}}_{\BS},A) \to H^m(\Gamma \bs D_{\Sig}, A)$$ 
for $m=0$ (resp.\ $m=0,1$, resp.\ $m\in \bZ$).
\end{cor}

\bigskip

\noindent {\rm Kazuya KATO
\\
Department of mathematics
\\
University of Chicago
\\
Chicago, Illinois, 60637, USA}
\\
{\tt kkato@math.uchicago.edu}

\bigskip

\noindent {\rm Chikara NAKAYAMA
\\
Department of Economics 
\\
Hitotsubashi University 
\\
2-1 Naka, Kunitachi, Tokyo 186-8601, Japan}
\\
{\tt c.nakayama@r.hit-u.ac.jp}

\bigskip

\noindent
{\rm Sampei USUI
\\
Graduate School of Science
\\
Osaka University
\\
Toyonaka, Osaka, 560-0043, Japan}
\\
{\tt usui@math.sci.osaka-u.ac.jp}
\end{document}